\documentclass{article}
\usepackage{amsmath}
\usepackage{graphicx}
\usepackage{algorithm}
\usepackage{amsmath}
\usepackage{algpseudocode}
\usepackage{amsmath}  % AMS math package
\usepackage{amssymb}  % AMS symbol package
\usepackage{booktabs}
\usepackage{comment}
\DeclareMathOperator{\Uncertainty}{Uncertainty}
\DeclareMathOperator{\Margin}{Margin}

\begin{document}

\title{The New Era of Dynamic Pricing: Synergizing Supervised Learning and Quadratic Programming.}
\author{Gustavo Bramao, Ilia Tarygin}
\date{}
\maketitle

\begin{abstract}
    In this paper, we explore a novel combination of supervised learning and quadratic programming to refine dynamic pricing models in the car rental industry. We utilize dynamic modeling of price elasticity, informed by ordinary least squares (OLS) metrics such as p-values, homoscedasticity, error normality. These metrics, when their underlying assumptions hold, are integral in guiding a quadratic programming agent. The program is tasked with optimizing margin for a given finite set target. Our methodology includes developing a simulation environment to assess various pricing strategies and to benchmark against established price heuristics. While the components of our approach—supervised learning and quadratic programming—are individually well-known, their integrated application in the context of dynamic pricing is, to the best of our knowledge, an innovative venture. This strategy is aimed not only at enhancing productivity and profitability but also at reducing operational costs. Additionally, it offers a model that is both explanatory and adaptable, capable of maintaining its robustness amid market fluctuations.
\end{abstract}

\maketitle
\section*{I. Introduction}
Pricing strategy is a cornerstone for businesses across various sectors, profoundly influencing their success and market position. This strategy intricately balances multiple factors, including supply and demand dynamics, competitor pricing, brand positioning, perceived value, and overarching business strategies. Despite its critical importance, many companies still rely on traditional, manual approaches to pricing. These methods often depend on the intuition and experience of domain experts, supplemented to some extent by data-driven insights. However, a paradigm shift is emerging in this domain, led by more innovative companies. For instance, companies like Lyft have revolutionized their approach to pricing. By leveraging advanced reinforcement learning techniques, they have managed to automate their pricing policies effectively (Qin et al., 2022). This automation represents a significant departure from conventional practices, highlighting the potential of artificial intelligence (AI) in transforming core business processes. This paper delves into this transformative potential, specifically focusing on the integration of supervised learning with quadratic programming to enhance dynamic pricing models. Our exploration aims to bridge the gap between traditional, expertise-driven pricing strategies and data-driven, AI-powered approaches. We propose a novel framework that not only streamlines the pricing process but also adapts nimbly to market changes, thereby offering businesses a competitive edge in an increasingly data-centric world.

\maketitle
\section*{II. Defining the Simulation Environment}
In this section, we detail the construction of a robust simulation environment designed to evaluate diverse pricing strategies. This environment is pivotal for the refinement of dynamic pricing models to reflect real-world market conditions accurately.

\section{Overview of the Simulation Environment}

Our simulation environment is designed to capture the demand forecast for a specific pick-up date at a baseline price. This baseline price is defined as a heuristic price, represented by the mean price derived from a set of rules. Mathematically, let $P_b$ denote the baseline price, which can be expressed as:
\begin{equation}
    P_b = \frac{1}{N} \sum_{i=1}^{N} p_i,
\end{equation}
where $N$ is the number of rules and $p_i$ is the price associated with the $i^{th}$ rule. Typically, $P_b$ is influenced by factors such as the length of rent $L$, advanced booking time $T$, branch type $B$, car group $C$, and whether it is a pick ($D_p$) or non-pick day ($D_{np}$). These factors can be represented as:
\begin{equation}
    F = \{L, T, B, C, D_p, D_{np}\}.
\end{equation}

Each combination of these factors is evaluated, and defining the grouping is a non-trivial task that requires extensive feature engineering. Arguably, one of the most crucial aspects of the success of this framework is how we define these groupings. Our aim is to be granular enough to obtain detailed cost and margin information while maintaining a good generalization of elasticity, minimizing cross-elasticity, and achieving accurate demand forecasting. In our framework, we aim to optimize using uncertainties that arise from the elasticity and demand forecast of the baseline price. The optimization can be represented as:

We aim to optimize the following function to achieve a balance between maximizing margin and minimizing uncertainty:
\begin{equation}
    \min_{\{F\}} \; \alpha \cdot \Uncertainty(\{E\}, \{D\}) - \beta \cdot \Margin(\{M\}),
\end{equation}
where $\Uncertainty(\{E\}, \{D\})$ denotes the function measuring the uncertainties from elasticity $\{E\}$ and demand forecast $\{D\}$, and $\Margin(\{M\})$ represents the margin function dependent on the set of features $\{F\}$. The coefficients $\alpha$ and $\beta$ are weights that determine the relative importance of minimizing uncertainty versus maximizing margin.

Here, $F$ denotes the set of features influencing both the uncertainty and margin, $E$ represents elasticity, $D$ stands for demand forecast, and $M$ encapsulates factors affecting the margin.

From the simulation, we transition from forecasting demand for a pick-up date to analyzing "on-rents" - cars that are rented out on a given date. This process involves sequential decision making under uncertainty.

\vspace{1cm}

\section{Forecasting and Distribution Models}
In our simulation, we develop a model that combines two main factors to improve predictions for rental reservations as they approach a specific pickup date. These factors are the Length of Rent (LOR) and Advanced Booking Time (ABT). We focus our forecasting efforts on the critical 14-day period immediately before the pickup date, aiming for high accuracy in short-term predictions rather than broad, long-term forecasts.
For predicting these time-sensitive trends, we propose Long Short-Term Memory networks (LSTMs). LSTMs, pioneered by Hochreiter and Schmidhuber in 1997, are widely recognized for their effectiveness in time series forecasting. Our approach targets making accurate forecasts for a given pickup date. Specifically, we forecast the overall demand, including all offers, even those that do not lead to actual reservations. We then estimate the conversion rate (CvR) using the elasticity coefficient. This allows us to calculate the expected number of reservations across different offer distributions and covariates, based on a hypothetical 1 per cent price change, which serves as our standard baseline price for comparison.

\begin{align*}
    f_t &= \sigma(W_f \cdot [h_{t-1}, x_t] + b_f) & \text{(Forget Gate)} \\
    i_t &= \sigma(W_i \cdot [h_{t-1}, x_t] + b_i) & \text{(Input Gate)} \\
    \tilde{C}_t &= \tanh(W_C \cdot [h_{t-1}, x_t] + b_C) & \text{(Candidate Cell State)} \\
    C_t &= f_t * C_{t-1} + i_t * \tilde{C}_t & \text{(Cell State Update)} \\
    o_t &= \sigma(W_o \cdot [h_{t-1}, x_t] + b_o) & \text{(Output Gate)} \\
    h_t &= o_t * \tanh(C_t) & \text{(Output)}
\end{align*}

The output \( h_t \) of the LSTM model is primarily employed to forecast demand, utilizing a distinct set of input features \( x_t \). These input features \( x_t \) include all the factors described in section 2.1. The LSTM model is trained by fine-tuning its weights \( W_f, W_i, W_C, W_o \) and biases \( b_f, b_i, b_C, b_o \) with the objective of minimizing the forecasting errors. This process involves iteratively adjusting these parameters to better align the model's predictions with the historical reservation data. Through this training, the model learns to effectively interpret the nuanced relationships between the base price, LOR, ABT, car groups, branch types, peaks and non peaks, thereby enhancing its capability to provide precise short-term demand forecasts.

\section{Optimization Solver for Maximizing Margin with Minimum Risk}

Our problem can be formulated as a function with quadratic prorpreties, as such we decide to use an optimization method using a quadratic program. The optimization solver is designed to identify the optimal pricing strategy that maximizes margin while minimizing risk associated with uncertainties in demand forecasts and elasticity values. The quadratic program takes into account the available capacity and incorporates measures of uncertainty into the decision-making process.

\subsection*{Mathematical Formulation}

The objective is to maximize the margin, which is the difference between total revenue and total costs, considering the available capacity and the uncertainties in demand forecast and elasticity values.

\subsubsection*{Total Demand Calculation}

The total demand $D_{\text{total}}$ for each booking segment (long, medium, short) is calculated as follows:
\begin{equation}
D_{\text{total}} = \sum_{i \in \{\text{long, medium, short}\}} B_i \times (1 + E_i \times \Delta P_i)
\end{equation}
where:
\begin{itemize}
    \item $B_i$ represents the initial bookings for segment $i$.
    \item $E_i$ is the price elasticity of demand for each segment $i$.
    \item $\Delta P_i$ is the percentage change in price for segment $i$.
\end{itemize}

\subsubsection*{Total Revenue}

Total revenue $R$ is computed as the sum of the products of the new prices $P'_i$ and the new demand $D'_i$ for each segment:
\begin{equation}
R = \sum_{i \in \{\text{long, medium, short}\}} P'_i \times D'_i
\end{equation}

\subsubsection*{Total Cost}

The total cost $C$ is given by the product of the total new demand across all segments and the fixed cost per booking $C_{\text{booking}}$:
\begin{equation}
C = D_{\text{total}} \times C_{\text{booking}}
\end{equation}

\subsubsection*{Margin Calculation}

The margin $M$ is the difference between total revenue and total cost:
\begin{equation}
M = R - C
\end{equation}

\subsubsection*{Incorporating Capacity and Uncertainty}

To include capacity and uncertainty considerations, the optimization problem is formulated as a quadratic programming problem (Qiu et al., 2021) that seeks to find the set of price changes $\Delta P$ which maximizes the expected margin while accounting for capacity constraints and minimizing risk associated with uncertainty in demand and elasticity estimates.

\begin{equation}
\max_{\Delta P} \mathbb{E}[M] - \lambda \cdot \text{Var}[M] \quad \text{subject to} \quad D_{\text{total}} \leq \text{Capacity}, \quad \Delta P_{\min} \leq \Delta P \leq \Delta P_{\max}
\end{equation}

where:
\begin{itemize}
    \item $\mathbb{E}[M]$ is the expected margin.
    \item $\text{Var}[M]$ is the variance of the margin, representing risk due to uncertainty.
    \item $\lambda$ is a risk aversion parameter that balances between margin maximization and risk minimization.
    \item Capacity constraints ensure that the total demand does not exceed the available capacity.
\end{itemize}

The optimization aims to balance the maximization of expected margin against the minimization of risk due to uncertainties in demand forecasts and elasticities, within the bounds of available capacity and permissible price adjustments.

\vspace{1cm}

Example of an optimal policy solved with QP:

\includegraphics[width=1\textwidth]{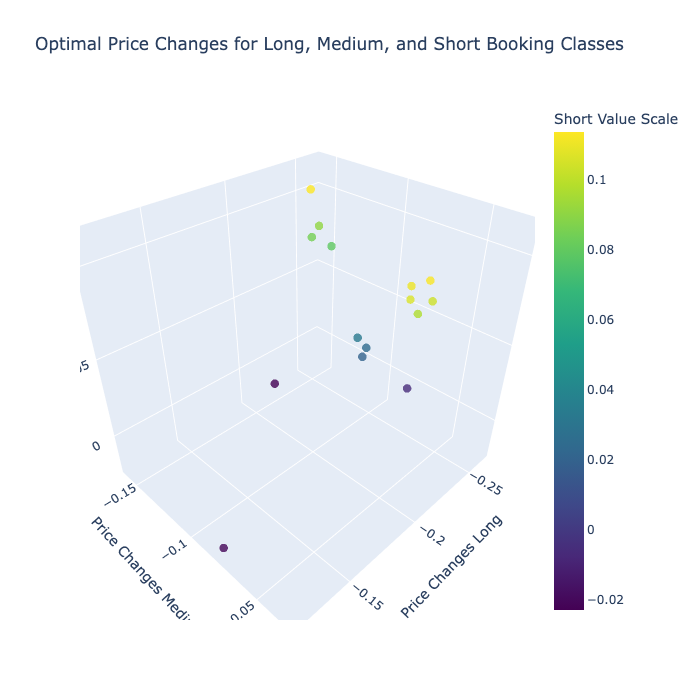}

\subsubsection*{Risk Formulation}

To conceptualize risk within our framework, we aggregate the series of sequential decisions emanating from the quadratic program, which is aimed at maximizing the cumulative margin. This approach inherently acknowledges the multifaceted nature of uncertainty within the decision-making process. Specifically, uncertainty originates from two primary sources: the standard deviation associated with the booking forecast at the heuristic price, denoted as \( \sigma_{BF} \), and the standard deviation related to the elasticity values predictions, \( \sigma_{E} \), obtained through our supervised learning model. A policy that is characterized by a higher degree of risk typically involves selecting a path of sequential decisions for booking classes that exhibit elevated uncertainties in both the booking forecast and the prediction of elasticities. The mathematical abstraction of risk, denoted as \( \mathcal{R} \), is formulated as follows:

\[
\mathcal{R} = \sum_{i} \sqrt{\sigma_{BF,i}^2 + \sigma_{E,i}^2},
\]

where:

\begin{itemize}
    \item \( \mathcal{R} \) represents the total risk associated with a sequence of decisions. It encapsulates the cumulative measure of uncertainty across all decisions made within the policy or strategy.
    \item \( \sigma_{BF,i} \) is the standard deviation of the booking forecast for decision \( i \) at the heuristic price, reflecting the variability and uncertainty associated with demand forecasts.
    \item \( \sigma_{E,i} \) represents the standard deviation of the elasticity values prediction for decision \( i \), capturing the uncertainty in the model's understanding of how demand responds to price changes.
\end{itemize}

\subsubsection*{Opportunity Cost}

Opportunity Cost is a measure used to evaluate the relative efficiency of different choices in decision-making. It is defined as the difference in value between the optimal choice and all other feasible alternatives. The optimal choice, in this context, is the one that maximizes margin, reduces risk under a certain capacity constraint thus, by definition, has an opportunity cost of zero since there is no better alternative forgone. Mathematically, for a set of possible choices \(\{C_1, C_2, \ldots, C_n\}\), with \(V(C_i)\) representing the value of choice \(i\), the opportunity cost of choosing any option \(C_k\) can be expressed as:

\[ OC(C_k) = \max\{V(C_1), V(C_2), \ldots, V(C_n)\} - V(C_k) \]

The choice \(C_m\) that satisfies \(\max\{V(C_1), V(C_2), \ldots, V(C_n)\} = V(C_m)\) is deemed optimal and its opportunity cost is \(OC(C_m) = 0\), indicating that it is the most efficient decision among all available options, with no value lost in the decision process.

\subsubsection*{Second Property of Opportunity Cost in Car Rental Operations}

A distinctive aspect of opportunity cost within the car rental industry lies in its operational flexibility, contrasting with the fixed capacities characteristic of airlines and hotels. This flexibility is exemplified by the potential to adjust fleet allocation in response to varying demand, such as relocating vehicles from downtown areas to airports.

To quantitatively assess this property, we propose a methodology that involves two computational steps. Initially, we execute a quadratic programming model subject to the existing fleet capacity constraints to establish a baseline. Subsequently, we run the same model without these constraints, simulating an ideal scenario where fleet size can be dynamically adjusted to perfectly match demand. The difference in profit margin between these two scenarios, normalized by the number of vehicles that could hypothetically be reallocated, quantifies the opportunity cost of not exploiting this operational flexibility.

Let \(P_c\) and \(P_u\) represent the profit margins obtained from the constrained and unconstrained quadratic programming models, respectively. If \(N\) denotes the number of vehicles that could potentially be reallocated, the opportunity cost per vehicle (\(OC_{\text{vehicle}}\)) can be mathematically defined as:

\begin{equation}
OC_{\text{vehicle}} = \frac{P_u - P_c}{\Delta N}
\end{equation}

where:
\begin{itemize}
    \item \(P_c\) = Profit margin under fixed capacity constraints,
    \item \(P_u\) = Profit margin with capacity unconstrained (ideal scenario),
    \item \(\Delta N = N_{\text{optimal}} - N_{\text{constrained}}\), the difference in the number of vehicles between the optimal and constrained scenarios.  
\end{itemize}

This formula provides a clear, quantifiable measure of the opportunity cost associated with not utilizing the car rental fleet's inherent flexibility to its fullest extent, offering insights into the economic implications of operational fleet management decisions.

\subsubsection*{Quadratic Programming Daily Batch}

Some of the interesting properties of quadratic programming (QP) include the existence of only one optimal solution and the computational efficiency in finding the optimal pricing policy. As such, every single time the program compiles, we can assume that the program will always find the optimal solution. We propose, therefore, running the QP as a batch every day to find the optimal policies. This will ensure that, regardless of the delta between predictions and observations, they will always be recalibrated based on actual observations. It is our belief that daily observations provide enough data for our context, limiting volatility and potentially overreacting to noise. 

\maketitle

\section*{III. Supervised Learning Methods for Addressing Endogeneity in Price-Demand Modeling}

In the car rental industry, estimating elasticity values for Length of Rental (LOR) and Advance Booking Time (ABT) faces the endogeneity challenge. Typically, increased demand causes agents to raise prices, mistakenly implying that higher prices boost demand. This misinterpretation is a classic example of the endogeneity problem.

To counter this, we implement a price randomization method, systematically altering prices independent of demand fluctuations. This strategy provides insights into the true impact of price changes on demand, eliminating the distortions caused by endogenous factors.

We then apply ordinary least squares (OLS) regression to analyze this randomized data. By thoroughly considering OLS's assumptions and inherent limitations, we ensure the reliability and accuracy of our elasticity calculations. These elasticity values are vital in our dynamic pricing model, especially within the quadratic programming framework, where they inform our pricing strategies. Our method enhances precision in understanding demand's response to price changes compared to traditional models.

\subsubsection*{Three-Step Approach:}

\textbf{Data Filtering:} \\
\textbf{Initial Data Set:} We start with an original dataset, represented as \(D\). \\
\textbf{Filtering Function:} A function \(F\) is defined, which evaluates each row \(x\) in \(D\). The function \(F(x)\) returns true if \(x\) meets specific pre-set criteria, such as defined ranges for LOR and ABT. \\
\textbf{Creation of Filtered Data Frame:} The filtered dataset, \(D_{\text{filtered}}\), consists of all rows \(x\) from \(D\) that satisfy \(F(x)\). 

\textbf{Log-Log Regression Modeling:} \\
After filtering, we apply a log-log regression model to \(D_{\text{filtered}}\) to estimate the interplay between price changes and demand. \\
The model is formulated as:
\begin{equation}
\log(Y) = \beta_0 + \beta_1\log(X) + \varepsilon
\end{equation}

\textbf{Causality and Strict Exogeneity:} \\
Under a specific set of narrow conditions, where the assumptions of Ordinary Least Squares (OLS) regression are met without violation, it becomes feasible to transition from correlation to causation. This capability stands as a prominent advantage of OLS methodology over alternative machine learning approaches (Schölkopf et al., 2021). From our standpoint, the primary objective is not to enhance out-of-sample predictions. Rather, our focus is on uncovering causal relationships that generalize well the effects of interventions. We aim to comprehend and predict the outcomes of price adjustments, specifically investigating the effects of both decreasing and increasing them, this involves measuring a decision under the influence of interventions, temporal shifts, and answering counterfactual questions through a causal lens. (Lewis and Wong, 2018).

\subsection*{Key Assumptions and Considerations in Our Approach:}

\subsubsection*{3.1 Addressing Endogeneity Bias:}
By incorporating random data, we reduce the chance of endogeneity bias. This bias arises when the variables we use to predict are correlated with the model's error term. In our situation, where changes in demand directly lead to price changes, this consideration is particularly relevant.

\subsubsection*{3.2 Compliance with OLS (Ordinary Least Squares) Assumptions:}
\textbf{Linearity:} The relationship between the independent and dependent variables should be linear. Our approach assumes that this linearity holds true in the log-log transformation of the variables. \\
\textbf{Independence:} Observations need to be independent of each other. Randomization helps in ensuring that the samples are independent and identically distributed (IID). \\
\textbf{Homoscedasticity:} Assumes that the error terms in our regression model have a constant variance across all observations. When this assumption is violated, resulting in heteroscedasticity (where error variances vary), it may lead to less efficient estimates and issues with statistical inference. To mitigate this, we aim to select predictors that demonstrate minimal heteroscedasticity, while also employing techniques such as robust standard errors to address potential variance inconsistencies in our model. \\
\textbf{Normal Distribution of Errors:} For the purposes of making inferences about the coefficients, the error terms should ideally follow a normal distribution.

\subsection*{3.3 Information gain versus entropy}

On the feature engineering aspect, we decide to take an approach that balances the causality of OLS and rigorousness of statistical inference, coupled with a hierarchical heuristic decision tree, that implicitly aims at balancing the information gain and minimizing the entropy. Our viewpoint on this is that "it is better to be overall right than precisely wrong". The logic is to go a subset more granular when we have enough evidence on the precision in the elasticity, we typically selected a p value less than 0.01 and an elasticity of variance of 1 point.

\subsection*{3.4.1 Predicting Future Elasticity: Bayesian Time Varying Coefficient Model}

We propose a method to forecast elasticity using a Bayesian time-varying coefficient that captures the trend to predict the future value of elasticity. Despite being simple in approach, we believe that elasticity isn't static in nature. Furthermore, we see value in accurately predicting global elasticity, or the highest node typically at the country level. When we lack sufficient information about a set of features, we fallback to the global elasticity. As such, it's our belief that predicting elasticity, rather than relying on past elasticity, will enhance the accuracy and efficiency of our quadratic programming, and increase margin by providing less risk and uncertainty to a set of features.

Given the single predictor in our model - the log of offer\_multiplier representing price changes from -15\% to +15\%, and the time component being the months, the BTVC model is formulated as follows:

\begin{equation}
    \beta_{t} = \text{function of time}(t), \quad t = 1, 2, \ldots, T
\end{equation}

where:
\begin{itemize}
\item $\beta_{t}$ represents the time-varying coefficient for the log of offer\_multiplier at week $t$.
\item $T$ is the total number of months observed.
\end{itemize}

\subsection*{Priors in the Bayesian Model}

In the Bayesian framework, we specify priors for the model parameters to incorporate existing knowledge and beliefs. For the intercept term $\alpha$ and the time-varying coefficients $\beta_t$, the following priors are considered:

\begin{itemize}
    \item Each time-varying coefficient \(\beta_t\) is modeled as following a normal distribution. This distribution is centered around the mean elasticity, denoted as \(\mu_{\text{{elasticity}}}\), and it has a variance \(\sigma^2_\beta\) that allows for flexibility in the estimates:
    \[\beta_t \sim \mathcal{N}(\mu_{\text{{elasticity}}}, \sigma^2_\beta), \quad \text{for } t = 1, 2, \ldots, T\]
    
\item The standard deviation of residuals, $\sigma$, is modeled with an Exponential distribution, a common choice for scale parameters in Bayesian models:
\[\sigma \sim \text{Exponential}(\lambda_{\sigma})\]
\end{itemize}

\subsection*{Implementation and Analysis using Stan}

The Bayesian Time Varying Coefficient (BTVC) model is implemented and analyzed using the Stan probabilistic programming language, interfaced through CmdStanPy in Python. This approach is particularly suited for handling the complexities inherent in Bayesian inference.

The core of the model lies in representing the log of reservations as a function of both time and the log-transformed offer multipliers. This relationship is encapsulated as follows:

\begin{align*}
    \text{for each time } t: \quad y_t &\sim \mathcal{N}(X_t \beta_t + \alpha, \sigma)
\end{align*}

In this formulation, \(y_t\) represents the log reservations for time \(t\), modeled as a normally distributed variable with a mean given by the linear combination of the predictors (log offer multipliers, \(X_t\)) and the time-varying coefficients (\(\beta_t\)), plus the intercept (\(\alpha\)). The standard deviation of the residuals is denoted by \(\sigma\).

The model is defined and compiled in Stan, leveraging Markov Chain Monte Carlo (MCMC) sampling to infer the posterior distributions of the model parameters, including the time-varying coefficient \(\beta_t\), inferring how the elasticity evolves.

\subsection*{Interpreting the Results}

The interpretation of the BTVC model primarily focuses on how the elasticity value varies over time. In the field of time-varying price elasticity, some research (Ruan et al., 2022), suggests simplistic models for time-varying coefficients by independently modeling $\Delta Q / \Delta P$ at each specific time, $t$. Our approach also incorporates this methodology for past observations.

For future predictions, we have chosen the Prophet library developed by Meta to project forthcoming trends, incorporating seasonality, trends, and cycles as primary components of time series analysis. Furthermore, we will revenue per day (RPD) as an additional regressor. Given our knowledge of anticipated daily revenue—which serves as both the outcome of our simulation and our objective—we aim to enhance our forecasting accuracy by integrating RPD into our model. To evaluate the effectiveness of our forecasts, we will benchmark them against two baselines: the mean of the observed data and the previous time point (t-1), under the assumption that t-1 provides a reliable indicator for future outcomes.

\begin{equation}
    y(t) = g(t) + s(t) + h(t) + \beta \cdot R(t) + \epsilon_t
\end{equation}

where:
\begin{itemize}
\item $y(t)$ represents the forecasted value of the target variable at time $t$.
\item $g(t)$ is the trend component at time $t$, modeling non-periodic changes.
\item $s(t)$ is the seasonality component at time $t$, capturing periodic changes.
\item $h(t)$ accounts for the effect of holidays or events at time $t$.
\item $\beta \cdot R(t)$ represents the impact of an external regressor, specifically daily revenue $R(t)$, on the forecast at time $t$, with $\beta$ being the coefficient measuring this impact.
\item $\epsilon_t$ is the error term at time $t$, capturing random fluctuations that are not explained by the model.
\end{itemize}

We employ a cross-validation method using historical data, with a 14-day horizon window, aiming for precise short-term elasticity forecasts. Accepting a trade-off for lower long-term performance, we use a rolling window technique. After each validation step, we extend the training period by 100 days for the next phase, testing the model's adaptability to new data and resilience to pattern changes over time.

\vspace{10pt}

\begin{tabular}{lc}
    \toprule
    Metric & Value \\
    \midrule
    RMSE for mean Static Elastic Value & 0.51 \\
    RMSE for Bayesian\_coefficients with Prophet & 0.18 \\
    RMSE for Naive Bayesian\_coefficients at t-1 & 0.21 \\
    \bottomrule
\end{tabular}

\vspace{10pt}

This indicates that with a limited amount of feature engineering applied to the Prophet model, and by leveraging the RPD as a regressor, it is possible to outperform the naive t-1 estimate. This finding also clearly illustrates that in situations of elasticity drift, utilizing the mean of the coefficient is not the most effective strategy. We believe that for predicting global elasticity, which acts as the highest fallback node when insufficient information is available, forecasting future values of elasticity can significantly increase margins by minimizing the cost associated with incorrect predictions. It is important to recognize that elasticity is arguably the most critical predictor in our quadratic programming model.

\section*{Conclusion}

In conclusion, our research introduces a modern approach by synergizing supervised learning and quadratic programming to refine dynamic pricing models, particularly in the car rental industry. To the best of our knowledge, this method represents a novel venture in dynamic pricing, blending the precision of supervised learning with the strategic optimization capabilities of quadratic programming. The core of our methodology lies in its pragmatic application and efficiency, it is our belief that when coupled with careful feature engineering, this method could significant increase margin. Our framework is distinguished by its explainability and interpretability, traits that are increasingly valued in the application of artificial intelligence in business contexts. By providing a clear rationale behind pricing decisions, stakeholders can gain insights into the model's operation, enhancing trust and facilitating strategic alignment. Moreover, our approach is computationally elegant, leveraging the strengths of both supervised learning for understanding demand elasticity and quadratic programming for optimizing pricing strategies under uncertainty.
A critical aspect of our method is its adeptness in managing risk. Through the careful consideration of uncertainties inherent in demand forecasts and elasticity estimates, our model smartly balances the objective of maximizing margins with the need to minimize risk. The introduction of the risk aversion parameter {$\lambda$} allows the decision maker to balance risk and margin, demonstrating a prudent approach to decision-making in volatile markets.

\section*{Appendix}
\subsection*{Granularity of the problem}
The key aspect to consider while defining an optimization model is granularity, it comes in the following dimensions:

\begin{itemize}
    \item 
    $N$ -- size of the pickup window, days;

    \item 
    $M$ -- maximum advanced booking time (ABT, time prior to pickup, antecedence, lead time), days;

    \item 
    $L$ -- lenght of rent (LoR), days;
    %\item 
    %$l \in \{\mbox{ECMR}, \mbox{CCMR}, ...\}$ -- car group 
\end{itemize}

\subsection*{Input variables}

\begin{itemize}
    \item 
    $F_i \in \mathbb{N}$ -- available fleet per pickup day;
    \item
    $P_{ijk}, C_{ijk}\in \mathbb{R^+}$ - average price and margin per demand group;    
    \item
    $D_{ijk} \in \mathbb{R^+}$ -- demand forecast at mean historical price as predicted number of reservations (per each day in the booking window, lenght of rent, ABT)
    \item
    $\epsilon_{jk} \in \mathbb{R^-}$ -- elasticity coefficients;
    \item
    $u_0\in \mathbb{R}_{[0, 100]}$ - expected utilization, percents.
\end{itemize}

\subsection*{Key equations}

Given the price relative change~$\Delta P$, the demand changes according to linear elasticity model:
\begin{equation}
    Demand(P\Delta P) = Demand(P)\left(1-elasticity(1-\Delta P)\right),
\end{equation}
hence margin at modified price for a specific demand slice can be expressed as:
\begin{eqnarray}
    Margin(P\Delta P) = Demand(P\Delta P)P\Delta P  - Costs =\\= Demand(P)\left(1-elasticity(1-\Delta P)\right)P\Delta P - Costs.
\end{eqnarray}
Transformation of reservations into utilization:
\begin{equation}
    u_t = \frac{1}{F_t}\sum_{i=0}^{t} \sum_{j=0}^{t-i} \sum_{k=t-i}^{L}D_{ijk}\left(1-\epsilon_{jk}\left(1 - \Delta P_{ijk}\right)\right), \,\, t \in \{1,...,N\},
\end{equation}
it's worth to mention that utilization is a linear function of price multipliers. 

\subsection*{Decision variables}
In an optimization problem, decision variables are the variables that can be adjusted or chosen by the solver in order to optimize a certain objective function. These variables represent the quantities or values that the decision-maker can control or influence in order to achieve the desired outcome. In our case they are~$\Delta P_{ijk}$ -- price multipliers.

\subsection*{Optimal pricing as a quadratic optimization problem}
The problem of finding an optimal pricing policy to maximize total margin in a given window of~$N$ booking days can be expressed as an quadratic optimization problem in the following way: 

\begin{flalign*}
    & \mathop{\mbox{maximize}}_{\Delta P} \sum_{i=0}^{N} \sum_{j=0}^{M} \sum_{k=0}^{L}D_{ijk}\left(P_{ijk}\left(\left(1-\epsilon_{jk}\right)\Delta P_{ijk}+\epsilon_{jk}\Delta P_{ijk}^2\right) - C_{ijk}\right)\\  
    & \mbox{subject to}\\
    & a u_0 \leq u_i \leq b u_0, \quad  i \in \{1,\dots,N\}  \\ 
    & \Delta P_{min} \leq \Delta P_{ijk} \leq \Delta P_{max}, \quad i,j,k \in \{1,\dots,N\}\times\{1,\dots,M\}\times\{1,\dots,L\}.
\end{flalign*}
Here~$a$ and~$b$, as well as~$\Delta P_{min}$ and~$\Delta P_{max}$ are some parameters to be defined heuristically or by the user.  
\section*{Dealing with uncertainty}
\begin{figure}[h!]    
    \includegraphics[width=0.9\textwidth]{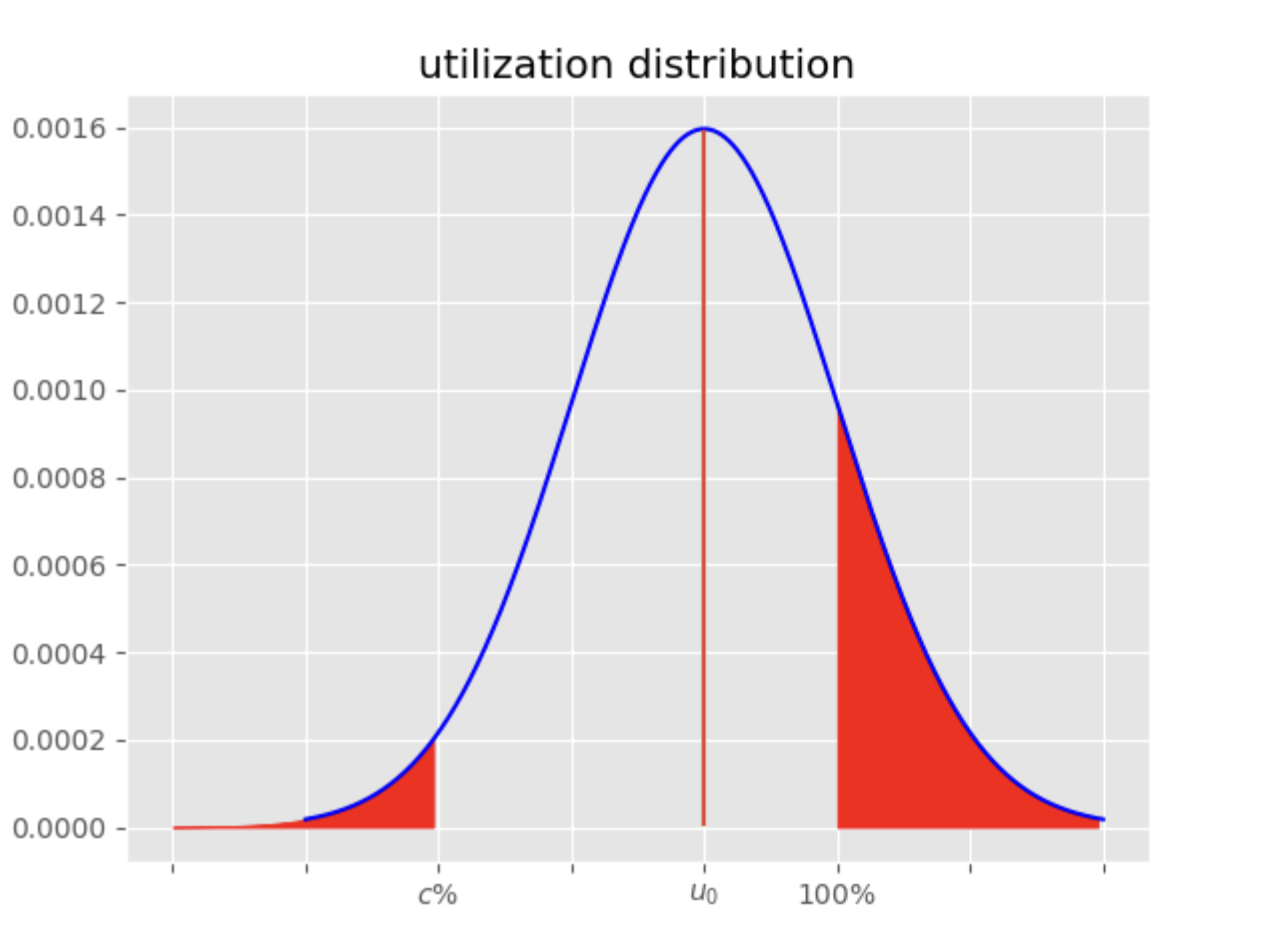}
    \caption{utilization risk as an area under the normal pdf.}
    \label{fig: utilization_risk}
\end{figure}

Let's define a \textit{risk} as a probability of some unpleasant event to happen, e.g.:
\begin{itemize}
    \item 
    probability of utilization on one of days in the pickup window being above or below some specific thresholds;
    \item 
    probability of total reservation margin to be below some specific threshold.
\end{itemize}
It's worth to mention, that demand forecast and elasticity are just estimations that are supplied with the measure of uncertainty. Let's make an assumption:
\begin{itemize}
    \item 
    $\widehat{D}_{ijk} \sim \mathcal{N}(D_{ijk},\,\sigma_{D_{ijk}}^{2})$;
    \item 
    $\widehat{\epsilon}_{ijk} \sim \mathcal{N}(\epsilon_{ijk},\,\sigma_{\epsilon_{ijk}}^{2})$,
\end{itemize}
then, given normally distributed independent variables have the properties
\begin{eqnarray}
    X \!&\sim& \!\mathcal{N}(\mu_{X},\,\sigma_{X}^{2}),\quad  Y \sim \mathcal{N}(\mu_{Y},\,\sigma_{Y}^{2}) \longrightarrow X+Y \!\sim \!\mathcal{N}(\mu_{X}+\mu_{Y},\,\sigma_{X}^{2}+\sigma_{Y}^{2})  \\ \nonumber
    X \!&\sim& \!\mathcal{N}(\mu_{X},\,\sigma_{X}^{2}),\quad  Y \sim \mathcal{N}(\mu_{Y},\,\sigma_{Y}^{2}) \longrightarrow X\cdot Y \!\sim \!\ p_Z(z) = \frac {K_0(|z|)}{\pi}, \;\;\; -\infty < z < +\infty,
\end{eqnarray}
where~$K_0$ refers to modified Bessel function of the second type,
we can state that utilization on a specific day is distributed normally with the standard deviation expressed as a function of price modifications:
\begin{equation}
    \sigma_{u_t}^2 = \frac{1}{F_t}\sum_{i=0}^{t}\displaylimits \sum_{j=0}^{t-i} \sum_{k=t-i}^{L}\left(1-\epsilon_{jk}\left(1 - \Delta P_{ijk}\right)\right)\sigma^2_{D_{ijk}},\quad t \in \{1,\dots, N\}.
\end{equation}
\begin{comment}
\begin{equation}
    M_{ijk} \!\sim \!\mathcal{N}(D_{ijk}\left(1-\epsilon_{ijk}\left(1-\Delta P_{ijk}\right)\right)P_{ijk}\Delta P_{ijk}, ).
\end{equation} 
\end{comment}
Thus the risk of the utilization on a specific date being above or below a specific threshold can be denoted as
\begin{equation}
    r_i = P\{u_i \leq c\} + P\{u_i > 100\} = 1 - \frac{1}{\sqrt{2\pi}}\int_{c}^{100}\displaylimits\!\!\exp({-\frac{u^2}{2\sigma_{u_i}^2}})du,
\end{equation}
where~$c$ is another parameter to be specied by the user. The illustration for the risk in given on \ref{fig: utilization_risk}.  

It is important that risk in our case comes as an explicit expression of decision variables, so it can be added to optimization problem as an additional constraint. For the normal cdf there is a well know approximation available:
\begin{equation*}
    \Phi(z) = \frac{1}{1+e^{-y}},
\end{equation*}
where~$y=1.526(1 + 0.1034z)$. Finally the risk of utilization being below a specified threshold~$c$ or above 1 can expressed as:
\begin{equation*}
    \Phi(z) \leq p,
\end{equation*}
where~$p$ is the measure of affordable risk.
In the same way we can express the risk of a total margin being below some specific threshold.

\end{document}